\documentclass[11pt,reqno,tbtags]{amsart}
\pdfoutput=1
\usepackage[]{amsmath,amssymb,amsfonts,latexsym,amsthm,enumerate}
\usepackage{amssymb,bbm,mathrsfs}
\usepackage{enumerate,graphicx,paralist}
\usepackage[numeric,initials,nobysame]{amsrefs}

\numberwithin{equation}{section}
\usepackage[colorlinks=true, pdfstartview=FitV, linkcolor=blue, citecolor=blue, urlcolor=blue,pagebackref=false]{hyperref}
\usepackage[marginratio=1:1,height=584pt,width=488pt,tmargin=117pt]{geometry}
\usepackage[margin=25pt,font=small,justification=centering
]{caption}

\newtheorem{maintheorem}{Theorem}

\def\P{\mathbb{P}}
\def\Z{\mathbb{Z}}
\def\R{\mathbb{R}}

\def\E{\mathbb{E}}

\newcommand{\cF}{\mathcal{F}}

\newcommand{\cW}{\mathcal{W}}

\begin{document}

\title{Note on the flux for TASEP with general disorder}

\author{Allan Sly}
\address{Allan Sly\hfill\break
Department of Mathematics\\
Fine Hall\\
Princeton University\\
Princeton, NJ 08540, USA.\hfill\break
Department of Statistics\\
UC Berkeley\\
Berkeley, CA 94720, USA.}
\email{asly@math.princeton.edu}
\urladdr{}

\begin{abstract}
Extending results of Bahadoran and Bodineau we show that the flux rate of TASEP with independent and identically distributed disorder always has a plateau of densities around $\frac12$.
\end{abstract}
\maketitle

The TASEP model (totally asymmetric simple exclusion process) on $\Z$ is a model of exclusion particles which always move to the right.  Configurations in $\{0,1\}^{\Z}$ represent the indicator of particles being present at each site on the line.  In TASEP with rate $r$ particles at site $i$ jump to the right according to a rate $r$ Poisson clock provided there is no particle present at site $i+1$.  A product measure with $\hbox{Bernoulli}(\rho)$ is stationary for all $\rho$ and the flux rate (the long term rate at which particles move across a bond) is given by
\[
f(\rho)  = r \rho(1-\rho).
\]
A natural problem is to consider the model where the jump rates of the bonds are chosen randomly from some law.  In this model at bond $i$ if a particle is present and $i+1$ is unoccupied then the particle jumps from $i$ to $i+1$  according to a rate $\alpha(i)$ Poisson clock for independent and identically distributed variables $\alpha(i)$.  This model was analysed in~\cite{BahBod16} in the case of a unit jump rate where an $\epsilon$ fraction of the sites have their jump rate replaced by some distribution $Q(\cdot)$ supported on $[r,1]$.
They show that the flux has a plateau around~$\frac12$, that is for some $0<\rho_c<\frac12$  it satisfies
\[
f(\rho) = r/4
\]
for all $\rho\in [\rho_c,1-\rho_c]$ provided $\epsilon$ is sufficiently small and $Q$ satisfies condition $\mathbf{H}$ which says that as $u \searrow 0$ then
\[
Q([r,r+u]) = O(u^\kappa),
\]
for some $\kappa > 1$.  They conjecture that this condition is not required and that the plateau should exist for any choice of $Q$ and in particular the most natural choice of $Q$ which is a point mass at $r$.  The proof makes use of a sophisticated multi-scale analysis.

In this note we show that the most general version of the statement can be proved significantly more simply.  Indeed if the jump rates $\alpha_i$ are independent and identically distributed and chosen from any positive non-constant distribution with infinimum $r>0$, i.e.
\[
r=\inf \{s: \P[\alpha \geq s] < 1\}
\]
then TASEP with site disorder $\alpha$ has a plateau.

\begin{maintheorem}
With disorder $\alpha$ and $\mu = \frac1{r} - \E[\frac1{\alpha(i)}]$ for all $\rho\in [\frac12 - \frac14 \mu r,\frac12 + \frac14 \mu r]$,
\[
f(\rho) = r/4.
\]
\end{maintheorem}

The simpler direction is the upper bound on the  flux of $r/4$ which follows from the existence of wide bands of sites with disorder close to $r$ as shown in Proposition~2.1 of~\cite{BahBod16}.

Now consider the lower bound on the flux.  The last passage percolation representation of TASEP is given as follows.  For  $\cW=\{(i,j )\in \Z^2:j\geq 0, i+j \geq 0\}$ let $Y_{ij}$ be independent rate $\alpha(i)$ exponential random variables.  A path $\gamma$ from $x$ to $y$ in $\cW$ is a sequence of vertices $x=(i_0,j_0),\ldots,(i_k,j_k)=y$ in $\cW$ such that
\[
(i_\ell-i_{\ell-1},j_\ell-j_{\ell-1}) \in \{(1,0),(-1,1)\}.
\]
That is each step moves either east or northwest.  For $(i,j)\in \cW$ we write $\Gamma_{i,j}$ for the set of all paths from $(0,0)$ to $(i,j)$ and let $T(i,j)$ denote the weight of the maximal path from the origin to $(i,j)$,
\[
T(i,j) = \max_{\gamma\in \Gamma_{i,j}} \sum_{(i',j')\in \gamma} Y_{i'j'}.
\]
We set $\cW'=\{(x,y )\in \R^2:y\geq 0, x+y \geq 0\}$.  For $(x,y)\in \cW'$ set
\begin{align*}
\tau(x,y)=  \lim_{n\to\infty} \frac1{n} T(\lfloor xn \rfloor, \lfloor yn \rfloor).
\end{align*}
In~\cite{Seppalainen99} it is shown that this limit exists almost surely and is deterministic (and in particular does not depend on the disorder $\alpha$).  It is concave, superadditive and positively 1-homogeneous.  We write
\[
h(t,x)=\inf\{y \geq 0 : \tau(x,y)> t\}
\]
which by homogeneity can be written as $h(t,x) = t k(x/t)$.  In the case of the the homogeneous TASEP with constant rate $r$ we have
\[
\tau(x,y)=\frac1{r}(\sqrt{x+y}+\sqrt{y})^2,\qquad k(v) = \begin{cases}
  \frac{r(1-v/r)^2}{4} &v\in[-1,1]\\
  -v &v \leq -1
\end{cases}
\]
It is shown in~\cite{Seppalainen99} that the flux function satisfies
\[
f(\rho) = \inf_{v\in \R} [k(v)+v\rho].
\]
To analyse $\tau$ and $h$ we will couple the process with a homogeneous rate $r$ TASEP.  Let $U_{ij}$ be independent random variables, independent of the $Y_{ij}$, with law $\hbox{Ber}(1-\frac{r}{\alpha(i)})\cdot \exp(r)$, that is the product of a Bernoulli and an exponential random variable.  Let $Z_{ij}=Y_{ij}+U_{ij}$ which is an exponential with rate $r$.  Thus $\E Y_{ij} = \frac1{r} - \E[\frac1{\alpha(i)}] = \mu>0$.
Let $\gamma_\star$ be the maximal path in the disorder system and so
\begin{align*}
\lim_{n\to\infty}\frac1{n} \E \max_{\gamma\in \Gamma_{i,j}}  \sum_{(i',j')\in \gamma} Y_{i'j'}
&= \lim_{n\to\infty} \frac1{n} \E \max_{\gamma\in \Gamma_{i,j}}  \sum_{(i',j')\in \gamma} Z_{i'j'} - \sum_{(i',j')\in \gamma} U_{i'j'}\\
&\leq \frac{1}{r} (\sqrt{x+y}+\sqrt{y})^2 - \liminf_n \frac1{n} \E[\sum_{(i',j')\in \gamma_\star} U_{i'j'}].
\end{align*}
Let $\cF$ be the $\sigma$-algebra generated by the disorder $(\alpha(i))$ and by $Y_{ij}$.  Since the $U_{ij}$ are conditionally independent given $\cF$ almost surely we have that
\begin{align*}
\liminf_n \frac1{n} \E[\sum_{(i',j')\in \gamma_\star} U_{i'j'}] &=  \liminf_n  \frac1{n}  \E[\E[\sum_{(i',j')\in \gamma_\star} U_{i'j'} \mid \cF]]\\
&= \liminf_n  \frac1{n}  \E[\sum_{(i',j')\in \gamma_\star} \E[ U_{i'j'} \mid \cF]]\\
& \geq \liminf_n  \frac1{n}  \E[\sum_{i=0}^{\lfloor n x \rfloor} \E[ U_{i0} \mid \cF]] =\mu |x|.
\end{align*}
where the second equality uses the fact that $\gamma_\star$ is $\cF$ measurable.  The inequality uses the fact that the $\E[ U_{i'j'} \mid \cF]$ are positive and depend only on $\alpha(i)$ and that there is at least one $(i',j')\in\gamma_\star$ for each $i' \in \{0,\ldots, \lfloor n x \rfloor\}$.  The final equality is the law of large numbers applied to $\E[ U_{i0} \mid \cF]$ which are bounded and independent and identically distributed.  Thus writing
\[
\tilde{\tau}(x,y)=\frac1{r}(\sqrt{x+y}+\sqrt{y})^2  - \mu |x|
\]
we have that
\[
\tau(x,y) \leq \tilde{\tau}(x,y).
\]
Now checking that
\[
f(\rho) = \inf_{v\in \R} [k(v)+v\rho]   \geq r/4
\]
is equivalent to
\begin{align}\label{e:tauStatement}
\inf_{v\in \R} [k(v)+v\rho]   \geq r/4
&\Longleftrightarrow \inf_{v\in \R} \ \frac{4}{r}k(v) \geq 1 - \frac{4}{r} v\rho \nonumber \\
&\Longleftrightarrow \inf_{x\in \R} \ \frac{4}{r}k(\frac{xr}{4}) \geq 1 - x\rho \nonumber \\
&\Longleftrightarrow \inf_{x\in \R} \ h(\frac{4}{r},x) \geq 1 - x\rho \nonumber \\
&\Longleftrightarrow \max_x \tau(x,1-x \rho) \leq \frac{4}{r},
\end{align}
where in the second line we made the substitution $v=\frac{xr}{4}$.  We will show that $\max_x \tilde{\tau}(x,1-x \rho) \leq \tilde{\tau}(0,1) = \frac{4}{r}$ for $\rho \in [\frac12 - \frac14 \mu r,\frac12 + \frac14 \mu r]$.  Note that the left and right derivatives of $\tilde{\tau}(x,1-x \rho)$ are respectively
\[
\lim_{h\nearrow 0}\frac{\tilde{\tau}(h,1-h \rho)-\tilde{\tau}(0,1)}{h} = \frac{2-4\rho}{r} +\mu,
\]
and
\[
\lim_{h\searrow 0}\frac{\tilde{\tau}(h,1-h \rho)-\tilde{\tau}(0,1)}{h} = \frac{2-4\rho}{r} -\mu.
\]
Thus the left derivative is non-negative and the right derivative is non-positive at $x=0$ when $\rho \in [\frac12 - \frac14 \mu r,\frac12 + \frac14 \mu r]$ then $\tau(x,1-x \rho)$ is maximized at $x=0$ since it is concave establishing equation~\eqref{e:tauStatement}.  Thus
\[
f(\rho) = \frac{r}{4} \quad \hbox{for } \rho \in [\frac12 - \frac14 \mu r,\frac12 + \frac14 \mu r].
\]

\bibliography{DisBib}
\bibliographystyle{plain}

\end{document}